\documentclass [a4paper,10pt]{amsart}
\usepackage{longtable}
\usepackage{hyperref}
\usepackage[T1]{fontenc}
\usepackage[utf8]{inputenc}
\usepackage[english]{babel}
\usepackage{textcomp}
\usepackage{dsfont}
\usepackage{latexsym}
\usepackage{amssymb}
\usepackage{amsthm}
\usepackage{amsmath}
\DeclareMathAlphabet{\mathpzc}{OT1}{pzc}{m}{en}
\usepackage{yfonts}
\usepackage{xfrac}
\usepackage{newlfont}
\usepackage{graphicx}
\usepackage{mathtools}
\usepackage{comment}
\usepackage{indentfirst}
\usepackage{braket}
\usepackage{mathrsfs}
\usepackage{xcolor}
\usepackage{fixmath}

\usepackage{etoolbox}

\usepackage{scalerel}[2014/03/10]
\usepackage[usestackEOL]{stackengine}
\newcommand{\dashint}{\,\ThisStyle{\ensurestackMath{%
			\stackinset{c}{.2\LMpt}{c}{.5\LMpt}{\SavedStyle-}{\SavedStyle\phantom{\int}}}%
		\setbox0=\hbox{$\SavedStyle\int\,$}\kern-\wd0}\int}

\DeclareMathOperator{\Hol}{Hol}

\DeclareMathOperator{\Int}{Int}

\renewcommand{\Re}{\mathrm{Re}\,}
\renewcommand{\Im}{\mathrm{Im}\,}

\newcommand{\vect}[1]{\mathbf{{#1}}}
\newcommand{\dd}{\mathrm{d}}

\DeclarePairedDelimiter{\abs}{\lvert}{\rvert}

\DeclarePairedDelimiter{\norm}{\lVert}{\rVert}

\let\originalleft\left
\let\originalright\right
\renewcommand{\left}{\mathopen{}\mathclose\bgroup\originalleft}
\renewcommand{\right}{\aftergroup\egroup\originalright}

\newcommand{\N}{\mathds{N}}

\newcommand{\C}{\mathds{C}}

\newcommand{\R}{\mathds{R}}

\newcommand{\T}{\mathds{T}}
\newcommand{\K}{\mathds{K}}

\newcommand{\cM}{\mathcal{M}}
\newcommand{\Nc}{\mathcal{N}}

\newcommand{\meg}{\leqslant}
\newcommand{\Meg}{\geqslant}
\newcommand{\eps}{\varepsilon}
\renewcommand{\phi}{\varphi}

\title{A Note on Hardy Spaces on Quadratic CR Manifolds}

\date{}
\begin{document}

\theoremstyle{definition}
\newtheorem{deff}{Definition}

\newtheorem{oss}[deff]{Remark}

\newtheorem{ass}[deff]{Assumptions}

\newtheorem{nott}[deff]{Notation}

\newtheorem{ex}[deff]{Example}

\theoremstyle{plain}
\newtheorem{teo}[deff]{Theorem}

\newtheorem{lem}[deff]{Lemma}

\newtheorem{prop}[deff]{Proposition}

\newtheorem{cor}[deff]{Corollary}

\author[M. Calzi]{Mattia Calzi}

\address{Dipartimento di Matematica, Universit\`a degli Studi di
	Milano, Via C. Saldini 50, 20133 Milano, Italy}
\email{{\tt mattia.calzi@unimi.it}}

\thanks{{\em Math Subject Classification 2020:} Primary: 32V20, 32A35; Secondary: 32A10,  32M10. }
\thanks{The author is a member of the	Gruppo Nazionale per l'Analisi Matematica, la Probabilit\`a e le	loro	Applicazioni (GNAMPA) of the Istituto Nazionale di Alta Matematica	(INdAM) and is  partially supported by the 2020	GNAMPA grant {\em Fractional Laplacians and subLaplacians on Lie groups and trees}.}

\begin{abstract} 
	Given a quadratic CR manifold $\cM$ embedded in a complex space, and a holomorphic function $f$ on a tubular neighbourhood of $\cM$, we show that the $L^p$-norms of the restriction of $f$ to the translates of $\cM$ is decreasing for the ordering induced by the closed convex envelope of the image of the Levi form of $\cM$. 
\end{abstract}
\maketitle

\section{Introduction}

Let $f$ be a holomorphic function on the upper half-plane $\C_+=\R+i \R_+^*$. If  $f$ belongs to the Hardy space $H^p(\C_+)$, that is, if $\sup_{y>0}\norm{f_y}_{L^p(\R)}$ is finite, where $f_y\colon x \mapsto f(x+iy)$, then it is well known that the function $y\mapsto \norm{f_y}_{L^p(\R)}$ is decreasing on $\R_+^*$, for every $p\in ]0,\infty]$. 
Nontheless, if $f$ is simply holomorphic, then the lower semicontinuous function $y\mapsto \norm{f_y}_{L^p(\R)}$ need not be decreasing. Actually, the set where it is finite may be any interval in $\R_+^*$, or even a disconnected set. 

Now, replace the upper half-plane $\C_+$ with a Siegel upper half-space 
\[
D\coloneqq \Set{(\zeta,z)\in \C^n\times \C\colon \Im z-\abs{\zeta}^2>0},
\]
and define
\[
f_h\colon \C^n\times \R \ni (\zeta,x)\mapsto f(\zeta,x+i\abs{\zeta}^2+h)
\]
for every $h>0$ and for every function on $D$. This definition is motivated by the fact that 
\[
b D\coloneqq \Set{(\zeta,x+i\abs{\zeta}^2)\colon (\zeta,x)\in \C^n\times \R}
\]
is the boundary of $D$, and the sets $b D+(0,i h)$, for $h>0$, foliate $D$ as the sets $\R+i y$, for $y>0$, foliate $\C_+$. 
If $f$ is holomorphic on $D$, then the mapping $h\mapsto \norm{f_h}_{L^p(\C^n\times \R)}$ is always decreasing  (though not necessarily finite), in contrast to the preceding case (cf.~Theorem~\ref{teo:1}).  This fact is closely related with the fact that evey holomorphic function defined in a neighbourhood of $b D$ automatically extends to $D$. More precisely, if one observes that $b D$ has the structure of a CR submanifold of $\C^n\times \C$, one may actually prove that every CR function (of class $C^1$) is the boundary values of a unique holomorphic function on $D$ (cf.~\cite[Theorem 1 of Section 15.3]{Boggess}). 

In this note we show that an analogous property holds when $b D$ is replaced by a general quadratic, or quadric, CR submanifold of a complex space, and then discuss some examples of \v Silov boundaries of (homogeneous) Siegel domains.

\section{Preliminaries}

We fix a complex hilbertian space $E$ of dimension $n$, a real hilbertian space $F$ of dimension $m$, and a hermitian map $\Phi\colon E\times E\to F_\C$. 
Define
\[
\cM\coloneqq \Set{(\zeta,x+i\Phi(\zeta))\colon \zeta\in E, x\in F }=\Set{(\zeta,z)\in E\times F_\C\colon \Im z-\Phi(\zeta)=0},
\]
where $F_\C$ denotes the complexification of $F$, while $\Phi(\zeta)\coloneqq\Phi(\zeta,\zeta)$ for every $\zeta\in E$. We define 
\[
\rho\colon E\times F_\C\ni (\zeta,z)\mapsto \Im z-\Phi(\zeta)\in F.
\]

We endow $E\times F_\C$ with the product
\[
(\zeta,z)(\zeta',z')\coloneqq (\zeta+\zeta', z+z'+2i \Phi(\zeta',\zeta))
\]
for every $(\zeta,z),(\zeta',z')\in E\times F_\C$, so that $E\times F_\C$ becomes a $2$-step nilpotent Lie group, and $\cM$ a closed subgroup of $E\times F_\C$. In particular, the identity of $E\times F_\C$ is $(0,0)$ and $(\zeta,z)^{-1}= (-\zeta,-z+2 i \Phi(\zeta))$ for every $(\zeta,z)\in E\times F_\C$.
It will be convenient to identify $\cM$ with the $2$-step nilpotent Lie group $\Nc\coloneqq E\times F$, endowed with the  product
\[
(\zeta,x)(\zeta',x')\coloneqq (\zeta+\zeta', x+x'+2 \Im \Phi(\zeta,\zeta'))
\]
for every $(\zeta,x),(\zeta',x')\in \Nc$, by means of the isomorphism
\[
\iota\colon\Nc \ni (\zeta,x)\mapsto (\zeta, x+i \Phi(\zeta))\in E\times F_\C.
\]
In particular, the identity of $\Nc$ is $(0,0)$ and $(\zeta,x)^{-1}=(-\zeta,-x)$ for every $(\zeta,x)\in \Nc$.
Notice that, in this way, $\Nc$ acts holomorphically (on the left) on $E\times F_\C$.
Given a function $f$ on $E\times F_\C$, we shall define
\[
f_h\colon \Nc\ni (\zeta,x)\mapsto f(\zeta,x+i \Phi(\zeta)+i h)\in \C
\]
for every $h\in F$.

Observe that the preceding groups structures show that, if we define the complex tangent space of $\cM$ at $(\zeta,z)$ as
\[
H_{(\zeta,z)}\cM\coloneqq T_{(\zeta,z)}\cM \cap (i T_{(\zeta,z)}\cM)
\]
for every $(\zeta,z)\in \cM$, where  $T_{(\zeta,z)}\cM$ denotes the real tangent space to $\cM$ at $(\zeta,z)$, identified with a subspace of $E\times F_\C$,
then
\[
H_{(\zeta,z)}\cM= \dd L_{(\zeta,z)} H_{(0,0)}\cM,
\]
where $L_{(\zeta,z)}$ denotes the left translation by $(\zeta,z)$ (in $E\times F_\C$), and $\dd L_{(\zeta,z)}$ its differential at $(0,0)$. Therefore, $\dim_\C H_{(\zeta,z)}=n$ for every $(\zeta,z)\in \cM$, so that $\cM$ is a CR submanifold of $E\times F_\C$ (cf.~\cite[Chapter 7]{Boggess}), called a qudratic or quadric CR manifold (cf.~\cite[Section 7.3]{Boggess} and~\cite{PelosoRicci,PelosoRicci2}).

We observe explicitly that $\cM$ is generic (that is, $\dim_\R \cM-\dim_\R H_{(0,0)}\cM=\dim_\R E\times F_\C-\dim_\R \cM$, cf.~\cite[Definition 5 and Lemma 4 of Section 7.1]{Boggess}) and that its Levi form may be canonically identified with $\Phi$ (cf.~\cite[Chapter 10]{Boggess} and~\cite{PelosoRicci2}).

\section{A Property of Hardy Spaces}

We denote by $C$ the convex envelope of $\Phi(E)$.

\begin{teo}\label{teo:1}
	Let $\Omega$ be an open subset of $F$ such that $\Omega=\Omega+\overline C$, and set $D\coloneqq \rho^{-1}(\Omega)$.
	Then, for every $f\in \Hol(D)$, for every $p\in ]0,\infty]$, for every $h\in \Omega$ and for every $h'\in \overline C$,
	\[
	\norm{f_{h+h'}}_{L^p(\Nc)}\meg \norm{f_h}_{L^p(\Nc)}.
	\]
\end{teo}

The proof is based on the `anaytic disc technique' presented in~\cite[Section 15.3]{Boggess}. 

Observe that the assumption that $\Omega=\Omega+\overline C$ is not restrictive. Indeed, if $\Omega$ is connected and $C$ has a non-empty interior $\Int C$, then every function which is holomorphic on $\rho^{-1}(\Omega)$ extends (uniquely) to a holomorphic function on $\rho^{-1}(\Omega+(\Int C\cup\Set{0}))$ by~\cite[Theorem 1 of Section 15.3]{Boggess}, and $\Omega+(\Int C\cup\Set{0})=\Omega+\overline C$ since $\Omega$ is open and $\overline C=\overline{\Int C}$ by convexity. The case in which $\Int C=\emptyset$ may be treated directly using similar techniques.

We also mention that, if $p<\infty$ and either  $\Phi$ is degenerate or the polar of $\Phi(E)$ has an empty interior (that is, the closed convex envelope of $\Phi(E)$ contains a non-trivial vector subspace), then either $f_h=0$ or $f_h\not \in L^p(\Nc)$ (at least for $p\Meg 1$ when $\Phi$ is non-degenerate). Cf.~\cite{Calzi} for more details in a similar case.

\begin{proof}
	For every $\vect v=(v_j)\in E^{m}$, consider 
	\[
	A_{\vect v}\colon \C \ni w \mapsto \bigg(\sum_{j=1}^{m} v_j w^j, i \sum_{j=1}^{m} \Phi(v_j)+ 2 i \sum_{k<j} \Phi(v_j, v_k) w^{j-k}\bigg)\in E\times F_\C, 
	\]
	and
	\[
	\Psi(\vect v)\coloneqq \sum_{j=1}^{m} \Phi(v_j)\in C,
	\]
	and observe that the following hold:
	\begin{itemize}
		\item $A_{\vect v}(0)= (0, i \Psi(\vect v))$;
		
		\item $\Psi(E^m)$ is the convex envelope of $\Phi(E)$, thanks to~\cite[Corollary 17.1.2]{Rockafellar};
		
		\item $\rho (A_{\vect v}(w))=0$ for every $w\in \T$;
		
		\item the mapping $A\colon E^{m}\ni\vect v\mapsto A_{\vect v}\in \Hol(\C; E\times F_\C)$ is continuous (actually, polynomial).
	\end{itemize}
	Now, take $h\in \Omega$. By continuity, there is $\eps>0$ such that $A_{\vect v}(\overline U)+i h\subseteq D$ for every $\vect v\in B_{E^{m}}(0,\eps)$, where $U$ denotes the unit disc in $\C$, and $\overline U$ its closure. Then, $A_{\vect v}(\overline U)+i h'\subseteq D$ for every $\vect v\in B_{E^{m}}(0,\eps)$ and for every $h'\in h+\overline C$.
	For every $h'\in \Psi( B_{E^{m}}(0,\eps))$, denote by $\nu_{h'}$ the image of the normalized Haar measure on $\T$ under the mapping $\pi\circ A_{\vect v}$, for some $\vect v \in B_{E^{m}}(0,\eps)\cap \Psi^{-1}(h')$, where $\pi\colon E\times F_\C\ni (\zeta,z)\mapsto (\zeta,x)\in \Nc$. Observe that, for every $(\zeta,x)\in \Nc$ and for every $h''\in h+\overline C$, the mapping
	\[
	\overline U \ni w \mapsto f((\zeta,x+i \Phi(\zeta))\cdot[ A_{\vect v}(w)+(0,i h'')])\in \C
	\]
	is continuous and holomorphic on $U$, so that, by subharmonicity (cf., e.g.,~\cite[Theorem 15.19]{Rudin}),
	\[
	\begin{split}
		\abs{f(\zeta,x+i \Phi(\zeta)+i(h'+h''))}^{\min(1,p)}&\meg \int_\T \abs{f((\zeta,x+i \Phi(\zeta))\cdot[ A_{\vect v}(w)+(0,i h'')])}^{\min(1,p)}\,\dd w\\
		&= \int_\Nc \abs{f_{h''}((\zeta,x)(\zeta',x'))}^{\min(1,p)}\,\dd \nu_{h'}(\zeta',x')\\
		&= \abs{f_{h''}}^{\min(1,p)}*\check\nu_{h'},
	\end{split}
	\]
	where $\check \nu_{h'}$ denotes the reflection of $\nu_{h'}$, while $\vect v$ is a suitable element of $ B_{E^{m}}(0,\eps)\cap \Psi^{-1}(h')$.
	Since $\nu_{h'}$ is a probability measure, by Young's inequality (cf., e.g.,~\cite[Chapter III, \S\ 4, No.\ 4]{BourbakiInt2}) we then infer that
	\[
	\norm{f_{h'+h''}}_{L^p(\Nc)}=\norm{\abs{f_{h'+h''}}^{\min(1,p)}}_{L^{\max(1,p)}}^{1/\min(1,p)}\meg \norm{\abs{f_{h''}}^{\min(1,p)}}_{L^{\max(1,p)}(\Nc)}^{1/\min(1,p)}=\norm{f_{h''}}_{L^p(\Nc)}
	\]
	for every $h'\in \Psi(B_{E^m}(0, \eps))$ and for every $h''\in h+\overline C$. Since every element of $C$ may we written as a finite sum of elements of $\Psi(B_{E^m}(0, \eps))$, the arbitrariness of $h''$ shows that
	\[
	\norm{f_{h+h'}}_{L^p(\Nc)}\meg\norm{f_{h}}_{L^p(\Nc)}
	\]
	for every $h'\in C$, hence for every $h'\in \overline C$ by lower semi-continuity. The proof is complete.
\end{proof}

\begin{cor}\label{cor:1}
	Assume that $C$ has a non-empty interior $\Omega$, and set $D\coloneqq \rho^{-1}(\Omega)$. Then, for every $p\in ]0,\infty]$ and $f\in \Hol(D)$,
	\[
	\sup_{h\in \Omega} \norm{f_h}_{h\in L^p(\Nc)}=\liminf_{h\to 0, h \in \Omega} \norm{f_h}_{L^p(\Nc)}.
	\]
\end{cor}

In particular, if we define the Hardy space $H^p(D)$ as the set of $f\in \Hol(D)$ such that $\sup_{h\in \Omega} \norm{f_h}_{h\in L^p(\Nc)}$ is finite,  the preceding result states that  $H^p(D)$ may be equivalently defined  as the set of $f\in \Hol(D)$ such that $\liminf\limits_{h\to 0, h \in \Omega} \norm{f_h}_{L^p(\Nc)}$ is finite. 
This result should be compared with~\cite{Boggess2}, where the boundary values of the elements of $H^p(D)$ are characterized as the CR elements of $L^p(\Nc)$, for $p\in [1,\infty]$. In particular, Corollary~\ref{cor:1} could be deduced from the results of~\cite{Boggess2}, when $p\in [1,\infty]$, though at the expense of some further technicalities. 

This result extends~\cite[Corollary 1.43]{CalziPeloso}.

\section{Examples}

We shall now present some exmples of homogeneous Siegel domains $D=\rho^{-1}(\Omega)$ for which $\overline{\Omega}$ is the closed convex envelope of $\Phi(E)$, so that Corollary~\ref{cor:1} applies.

We recall that $D$ is said to be a Siegel domain if $\Omega$ is an open convex cone not containing affine lines,  $\Phi$ is non-degenerate, and $\Phi(E)\subseteq \overline \Omega$. In addition, $D$ is said to be homogeneous if the group of its biholomorphisms acts transitively on $D$.
It is known (cf., e.g.,~\cite[Proposition 1]{Calzi2}) that $D$ is homogeneous if and only if there is a triangular Lie subgroup $T_+$ of $GL(F)$ which acts simply transitively on $\Omega$, and for every $t\in T_+$ there is $g\in GL(E)$ such that $t\Phi= \Phi(g\times g)$. 

If $T'_+$ is another Lie subgroup of $GL(F)$ with the same properties as $T_+$, then $T_+$ and $T'_+$ are conjugated by an automorphism of $F$ preserving $\Omega$.
Thanks to this fact, we may use the results of~\cite{CalziPeloso} even if a different $T_+$ is  chosen. In particular, there is a surjective (open and) continuous homomorphism  of Lie groups
\[
\Delta\colon T_+\to (\R_+^*)^r
\]
for some $r\in \N$, called the rank of $\Omega$, so that
\[
\Delta^{\vect s}=\Delta_1^{s_1}\cdots \Delta_r^{s_r},
\]
$\vect s\in \C^r$, are the characters of $T_+$.  Once a base point $e_\Omega\in \Omega$  has been fixed, $\Delta^{\vect s}$ induces a function $\Delta_\Omega^{\vect s}$ on $\Omega$, setting $\Delta_\Omega^{\vect s}(t(e_\Omega))=\Delta^{\vect s}(t)$ for every $t\in T_+$.

Up to modify $\Delta$, we may then assume that the functions $\Delta_\Omega^{\vect s}$ are bounded on the bounded subsets of $\Omega$ if and only if $\Re \vect s \in \R_+^r$ (cf.~\cite[Lemma 2.34]{CalziPeloso}).
In particular, there is $\vect b\in \R_-^r$ such that $\Delta^{-\vect b}(t)= \abs{\det_\C g}^2$ for every $t\in T_+$ and for every $g\in GL(E)$ such that $t\Phi=\Phi(g\times g)$  (cf.~\cite[Lemma 2.9]{CalziPeloso}), and one may prove that $\vect b\in (\R_-^*)^r$ if and only if $\Phi(E)$ generates $F$ as a vector space, in which case $\Omega$ is the interior of the convex envelope of $\Phi(E)$ (cf.~\cite[Proposition 2.57 and its proof, and Corollary 2.58]{CalziPeloso}). Therefore, we are interested in finding examples of homogeneous Siegel domains for which $\vect b\in (\R_-^*)^r$.

Notice, in addition, that if $\vect b\not \in (\R_-^*)^r$, then $\Phi(E)$ is contained in a hyperplane, so that the interior of its convex envelope is empty.

The Siegel domain $D$ is said to be symmetric if it is homogeneous and  admits an involutive biholomorphism with a unique fixed point (equivalently, if for every $(\zeta,z)\in D$ there is an involutive biholomorphism of $D$ for which $(\zeta,z)$ is an isolated (or the unique) fixed point). The domain $D$ is said to be irreducible if it is not biholomorphic to the product of two non-trivial Siegel domains. 

It is well known that every symmetric Siegel domain is biholomorphic to a product of irreducible ones, and that the irreducible symmetric Siegel domains can be classified in four infinite families plus two exceptional domains (cf., e.g.,~\cite[\S\S\ 1, 2]{Arazy}).  
In particular, for an irreducible symmetric Siegel domain, either $\vect b=\vect 0$ (that is, $E=\Set{0}$, in which case $D$ is `of tube type'), or $\vect b\in (\R_-^*)^r$ (cf., e.g.,~\cite[Example 2.11]{CalziPeloso}). Hence, when $D$ is a symmetric Siegel domain,  $\overline\Omega$ is the closed convex envelope of $\Phi(E)$  if and only if none of the irreducible components of $D$ is of tube type. Note that these domains can be also characterized as those which do not admit any non-constant \emph{rational} inner functions, thanks to~\cite{KoranyiVagi}.

\medskip

We now present some examples of (homogeneous) Siegel domains.

\begin{ex}
Let $\K$ be either $\C$ or the division ring of the quaternions. In addition, fix $r,k,p\in \N$ with $p\meg r$, and define
\begin{itemize}
	\item $E$ as the space of $k\times r$ matrices over $\K$ whose $j$-th columns have zero entries for $j=p+1,\dots, r$;
	
	\item $F$ as the space of self-adjoint $r\times r$ matrices over $\K$;
	
	\item $\Omega$ as the cone of non-degenerate positive self-adjoint $r\times r$ matrices over $\K$;
	
	\item 
	\[
	\Phi\colon E\times E \ni (\zeta,\zeta')\mapsto  \frac 1 2 [(\zeta'^*\zeta+\zeta^*\zeta')+i (\zeta^*i \zeta'-\zeta'^*i \zeta) ]\in F_\C;
	\]
	
	\item $T_+$ as the group of upper triangular $r\times r$-matrices over $\K$ with strictly positive diagonal entries, acting on $\Omega$ (and $F$) by the formula $t\cdot h \coloneqq t h t^*$;
	
	\item $\Delta\colon T_+\ni t \mapsto (t_{1,1},\dots,t_{r,r})\in (\R_+^*)^r$.
\end{itemize}
Then, $\Omega$ is an irreducible symmetric cone\footnote{A cone is said to be homogeneous if the group of its linear automorphisms acts transitively on it. It is said to be symmetric if, in addition, it is self-dual for some scalar product. A convex cone is said to be irreducible if it is not isomorphic to a product of non-trivial convex cones.} of rank $r$ on which $T_+$ acts simply transitively by~\cite[Example 2.6]{CalziPeloso}. In addition, $\Phi$ is well defined, since $\zeta'^*\zeta+\zeta^*\zeta', \zeta^*i \zeta'-\zeta'^*i \zeta\in F$ for every $\zeta,\zeta'\in E$,  and clearly $\Phi(\zeta)\in \overline \Omega$ and
\[
t\cdot \Phi(\zeta)= t\cdot (\zeta^*\zeta)= (\zeta t^*)^*(\zeta t^*)= \Phi( \zeta t^*)
\]
for every $t\in T_+$ and for every $\zeta\in E$ (with $\zeta t^*\in E$), so that $D$ is homogeneous. 
Then, $\vect b=(b_j)$, with $b_j=- k \dim_\C \K$ for $j=1,\dots,p$ and $b_j=0$ for $j=p+1,\dots, r$. 
Consequently, $\overline \Omega$ is the closed convex envelope of $\Phi(E)$  if and only if $p=r$ and $k>0$.

Notice that $D$ is irreducible since $\Omega$ is irreducible (cf.~\cite[Corollary 4.8]{Nakajima}), and that $D$ is symmetric if $k p=0$ or if $p=r$ and $\K=\C$ (cf.~\cite[Examples 2.14 and 2.15]{CalziPeloso}). If $k p (r-p)> 0$, or if $\K\neq \C$, $r\Meg 3$, and $k\Meg 2$, then $D$ cannot be symmetric.
\end{ex}

\begin{ex}
Take $k,p,q\in \N$, $p\meg 2$. Define:
\begin{itemize}
	\item $E$ as the space of formal $k\times 2$ matrices whose entries of the first column belong to $\C$ (and are $0$ if $p=0$), and whose entires of the second column belong to $\C^q$ (and are $0$ if $p\meg 1$);
	
	\item $F$ as the space of formally self-adjoint $2\times 2$ matrices whose diagonal entries belong to $\R$, and whose non-diagonal entries belong to $\C^q$;
	
	\item $\Omega$ as the cone of $\left( \begin{smallmatrix}
		a & b\\ \overline b & c
	\end{smallmatrix} \right)\in F$ with $a,c>0$, $b\in \C^q$, and $a c-\abs{b}^2>0$;

	\item $\Phi$ so that
	\[
	\Phi\left(\begin{matrix}
		a_1 & b_1 \\ \vdots & \vdots \\ a_k & b_k
	\end{matrix}  \right) = \left( \begin{matrix}
		\sum_j \abs{a_j}^2 & \sum_j \overline{a_j} b_j  \\
		\sum_j a_j \overline{b_j}& \sum_j \abs{b_j}^2
	\end{matrix}\right)
	\]
	for every $\left(\begin{smallmatrix}
		a_1 & b_1\\ \vdots &\vdots\\ a_k & b_k
	\end{smallmatrix}\right)\in E$;

	\item $T_+$ as the group of formal $2\times 2$ upper triangular matrices with diagonal entries in $\R_+^*$ and non-diagonal entries in $\C^q$, with the action\footnote{Formally, $\left( \begin{smallmatrix}
			a & b\\
			0 & c
		\end{smallmatrix}\right)\cdot \big( \begin{smallmatrix}
			a' & b'\\
			\overline b' & c'
		\end{smallmatrix}\big)=\left( \begin{smallmatrix}
		a & b\\
		0 & c
	\end{smallmatrix}\right) \big( \begin{smallmatrix}
	a' & b'\\
	\overline b' & c'
\end{smallmatrix}\big)\left( \begin{smallmatrix}
a & b\\
0 & c
\end{smallmatrix}\right)^*$.}
	\[
	\left( \begin{matrix}
		a & b\\
		0 & c
	\end{matrix}\right)\cdot \left( \begin{matrix}
	a' & b'\\
	\overline b' & c'
	\end{matrix}\right)\coloneqq \left( \begin{matrix}
	a' a^2+c'\abs{b}^2+2 a \Re \langle b,b'\rangle &  ac b'+cc' b\\
	a c \overline{b'}+c c' \overline b & c^2 c'
	\end{matrix}\right);
	\]  
	
	\item $\Delta\colon T_+\ni t \mapsto (t_{1,1}, t_{2,2})$.
\end{itemize}
Then, $\Omega$ is an irreducible symmetric cone of rank $2$ on which $T_+$ acts simply transitively (cf.~\cite[Example 2.7]{CalziPeloso}). In addition, $\Phi(\zeta)\in \overline \Omega$ for every $\zeta\in E$, and
\[
t\cdot \Phi(\zeta)=\Phi(\zeta t^*)
\]
for every $t\in T_+$ and $\zeta\in E$ (with $\zeta t^*\in E$), provided that $p\meg 1$. 
Then, $D$ is an irreducible Siegel domain, and it is homogeneous if $p\meg 1$ (it is symmetric if $p=0$). In addition, $\vect b=\vect 0$ if $p=0$, while $\vect b=(k,0 )$ if $p=1$.
Further, if $p=2$, then $\Phi(E)$ contains the boundary of $\Omega$, since  $\left(\begin{smallmatrix}
	a & b\\ 
	\overline b & c
\end{smallmatrix}\right)=\Phi\left(\begin{smallmatrix}
a^{1/2} &  a^{-1/2} b\\
0 &0\\
\vdots & \vdots\\
0 &0
\end{smallmatrix}\right)$, for every $a>0$, for every $c\Meg 0$ and for every $b\in \C^q$ such that $\abs{b}^2=a c$ (the case $a=0$, $b=0$, $c\Meg 0$ is treated similarly). Then, $\overline \Omega$ is the closed convex envelope of $\Phi(E)$ if and only if $p=2$.
\end{ex}


\begin{thebibliography}{9}
\bibitem{Arazy}
Arazy, J., \emph{A Survey of Invariant Hilbert Spaces of Analytic Functions on Bounded Symmetric Domains}, {Contemp.\ Math.} {185} (1995), p.~7--65.

	\bibitem{Boggess}
	Boggess, A., \emph{CR Manifolds and the Tangential Cauchy--Riemann Complex}, CRC Press, 1991.
	
	
	\bibitem{Boggess2}
	Boggess, A., CR Extension for $L^p$ CR Functions on a Quadric Submanifold of $\C^n$, \emph{Pac.\ J.\ Math.} \textbf{201} (2001), p.~1--18.
	
	\bibitem{BourbakiInt2}
	Bourbaki, N., \emph{Integration II}, Chap.\ 7--9, Springer, 2004.

	\bibitem{Calzi2}
	Calzi, M., Besov Spaces of Analytic Type: Interpolation, Convolution, Fourier Multipliers, Inclusions, preprint (2021), arXiv:2109.09402.
	
	\bibitem{Calzi}
	Calzi, M., Paley--Wiener--Schwartz Theorems on Quadratic CR Manifolds, preprint (2021), arXiv:2112.07991.
	
	\bibitem{CalziPeloso}
	Calzi, M., Peloso, M.\ M., Holomorphic Function Spaces on Homogeneous Siegel Domains, \emph{Diss.\ Math.} \textbf{563} (2021), p.~1--168.

	\bibitem{KoranyiVagi}
	Kor\'anyi, A., V\'agi, S., Rational Inner Functions on Bounded Symmetric Domains, \emph{Trans.\ Amer.\ Math.\ Soc.} \textbf{254} (1979), p.~179--193. 

	\bibitem{Nakajima}
	Nakajima, K., Some Studies on Siegel Domains, \emph{J.\ Math.\ Soc.\ Japan} \textbf{27} (1975), p.~54--75.
	
	\bibitem{PelosoRicci}
	Peloso, M.\ M., Ricci, F., Tangential Cauchy--Riemann Equations on Quadratic CR Manifolds, \emph{Rend.\ Mat.\ Acc.\ Lincei} \textbf{13} (2002), p.~125--134.
	
	\bibitem{PelosoRicci2}
	Peloso, M.\ M., Ricci, F., Analysis of the Kohn Laplacian on Quadratic CR Manifolds, \emph{J.\ Funct.\ Anal.} \textbf{203} (2003), p.~321--355.
	
	\bibitem{Rockafellar}
	Rockafellar, R.\ T., \emph{Convex Analysis}, Princeton University Press, 1970.
	
	\bibitem{Rudin}
	Rudin, W., \emph{Real and Complex Analysis}, McGraw-Hill, 1987.
\end{thebibliography}
\end{document}